\documentclass[12pt, reqno]{amsart}

\usepackage[utf8]{inputenc}
\usepackage{csquotes}
\usepackage[
backend=biber,
style=alphabetic,
doi=false,isbn=false,url=false]{biblatex}
\usepackage[english]{babel}
\usepackage[margin=3cm]{geometry}
\usepackage{amssymb}
\usepackage{amsthm}
\usepackage{amsmath}
\usepackage{mathtools}
\usepackage{marvosym}
\usepackage{bbold}
\usepackage[shortlabels]{enumitem}
\usepackage{amsaddr}
\usepackage{xcolor}
\usepackage{tikz-cd}
\usepackage[mathscr]{euscript}

\definecolor{blauD}{HTML}{233D4D}
\usepackage{hyperref}
\hypersetup{colorlinks=true, linkcolor=blauD, citecolor=blauD, urlcolor=blauD}
\usepackage[capitalize]{cleveref}

\newcommand{\Def}[1]{{\bf #1}}

\newtheorem{thm}{Theorem}[section]
\newtheorem{prop}[thm]{Proposition}
\newtheorem{lm}[thm]{Lemma}
\newtheorem{cor}[thm]{Corollary}
\theoremstyle{definition}

\newtheorem{smpl}[thm]{Example}
\newtheorem{rem}[thm]{Remark}

\crefname{lm}{Lemma}{Lemmas}
\crefname{thm}{Theorem}{Theorems}
\crefname{prop}{Proposition}{Propositions}
\crefname{defin}{Definition}{Definitions}
\crefname{rem}{Remark}{Remarks}
\crefname{smpl}{Example}{Examples}

\DeclareMathOperator{\im}{im}
\DeclareMathOperator{\rk}{rk}
\DeclareMathOperator{\cl}{cl}
\DeclareMathOperator{\WQSym}{\mathtt{WQSym}}
\DeclareMathOperator{\Mat}{\mathtt{Mat}}
\DeclareMathOperator{\ncPsi}{\text{\Neptune}}
\DeclareMathOperator{\Loops}{Loops}
\DeclareMathOperator{\Max}{Max}
\DeclareMathOperator{\type}{type}
\newcommand{\Mnco}{\text{\Scorpio}}
\newcommand{\bw}{\mathbf{w}}

\newcommand{\nd}{d}
\newcommand{\matroid}{\mathtt{M}}
\newcommand{\natroid}{\mathtt{N}}
\newcommand{\bases}{\mathcal{B}}
\newcommand{\unifMat}[2]{\mathtt{U}_{#1}^{#2}}
\newcommand{\opi}{\vec{\boldsymbol{\pi}}}
\newcommand{\otau}{\vec{\boldsymbol{\tau}}}
\newcommand{\flag}{\mathcal{C}}

\newcommand{\setofflags}{\mathfrak{F}}
\newcommand{\Q}{\mathbb{Q}}
\newcommand{\R}{\mathbb{R}}
\newcommand{\Z}{\mathbb{Z}}

\newcommand{\groundset}{U} 
\newcommand{\poset}{\mathscr{P}} 

\newcommand{\poly}{G}

\begin{document}

\title[Chromatic WQSym of nested matroids are independent]{Chromatic word-quasisymmetric functions of matroids II:
 Linear independence}

\author{Raul Penaguiao}
\email{raul.penaguiao@proton.me}
\address{University of Basel}
\author{Sophie Rehberg}
\email{rehberg.sophie@uqam.ca}
\address{Universit\'e du Qu\'ebec \`a Montr\'eal}

\date{November 26, 2026}

\begin{abstract}
In recent work (arXiv:2609.29927) the authors associated to every matroid $\matroid$ on $[\nd]$ a word-quasisymmetric function $\ncPsi(\matroid)$, showed that the image of $\ncPsi$ in degree $\nd$ is spanned by the images of the $\nd!$ loopless nested matroids, and conjectured that these images are linearly independent.
We prove this conjecture by appealing to a linear independence result from Hampe (2017).
This proof was developed using AI.
\end{abstract}

\maketitle

\section{Introduction}

For a matroid $\matroid$ on $[\nd]$ and a weight function $f\colon[\nd]\to\Z_{\ge1}$, the weight of a basis $B$ is $f(B)=\sum_{i\in B}f(i)$, and $f$ is called \Def{$\matroid$-generic} if exactly one basis of $\matroid$ has maximum weight.
The \Def{chromatic word-quasisymmetric function} of $\matroid$ is the formal power series in non-commuting variables $\bw_1,\bw_2,\dots$
\begin{equation*}
 \ncPsi(\matroid)\coloneqq\sum_{\substack{f\colon[\nd]\to\Z_{\ge1}\\ f \text{ is } \matroid\text{-generic}}} \bw_{f(1)}\bw_{f(2)}\cdots\bw_{f(\nd)}\,.
\end{equation*}
It was studied in \cite{penaguiao2026chromaticb}, as a non-commutative analogue of the quasisymmetric function of a matroid of Billera, Jia and Reiner \cite{billera2009quasisymmetric}.
Extending linearly, one obtains a linear map $\ncPsi\colon\Mat_\nd\to\WQSym_\nd$ from the vector space spanned by matroids on $[\nd]$ to word-quasisymmetric functions of degree $\nd$.
In \cite{penaguiao2026chromaticb} it is shown, using the valuativity of $\ncPsi$ and results of Derksen and Fink \cite{derksen2010valuative}, that the image $(\im\ncPsi)_\nd$ is spanned by $\ncPsi(\natroid)$ for $\natroid$ ranging over the $\nd!$ loopless nested matroids on $[\nd]$.
The main result of this note is that this spanning set is a basis.

\begin{thm}\label{thm:main}
Let $\nd\ge1$.
The family $\{\ncPsi(\natroid)\}_\natroid$, where $\natroid$ ranges over the loopless nested matroids on $[\nd]$, is linearly independent in $\WQSym_\nd$.
Consequently, $\dim(\im\ncPsi)_\nd=\nd!$.
\end{thm}

Together with \cite{penaguiao2026chromaticb}, this shows that the kernel of $\ncPsi$ in degree $\nd$ is spanned by the valuative relations together with the loop-coloop relations described there.

\subsection*{Idea of the proof}
Our main tool is a result by Hampe \cite[Proposition 3.2]{hampe_intersection_2017}: for a loopless matroid $\natroid$ of rank $r$ on $[\nd]$, record which chains $\flag\colon\emptyset\subsetneq F_1\subsetneq\cdots\subsetneq F_r=[\nd]$ consist of flats of $\natroid$; the resulting $\{0,1\}$-vectors are linearly independent as $\natroid$ ranges over the loopless nested matroids of rank $r$ (see also \cref{thm:hampe} below).

For every such chain $\flag$ of length $r$, we define a linear functional $\Lambda_\flag\colon\WQSym_\nd\to\Q$ with the following two properties (see \cref{prop:chainfunctional}): for a loopless matroid $\matroid$,
\begin{itemize}
 \item if $\rk(\matroid)>r$, then $\Lambda_\flag(\ncPsi(\matroid))=0$;
 \item if $\rk(\matroid)=r$, then $\Lambda_\flag(\ncPsi(\matroid))=1$ if all $F_i$ are flats of $\matroid$, and $0$ otherwise.
\end{itemize}
We order the loopless nested matroids by rank, and the functionals $\Lambda_\flag$ by the length $r$ of the chain.
More precisely, we choose some linear order for loopless nested matroids that respects the rank and some linear order for chains $\flag$ that respects the lengths of the chains.
Then the matrix $\big(\Lambda_\flag(\ncPsi(\natroid))\big)_{\natroid, \flag}$ with $\natroid$ ranging over loopless nested matroids and $\flag$ ranging over chains is block triangular, and its diagonal blocks consist of Hampe's chains of flats vectors, which are linearly independent.
This gives \cref{thm:main} which we will prove in \cref{sec:main}.
Note that nothing is claimed for $\rk(\matroid)<r$; indeed $\Lambda_\flag(\ncPsi(\matroid))$ can be nonzero in that case (\cref{smpl:lowrank}), and the triangular structure is all we need.

The functional $\Lambda_\flag$ is built from three ingredients.
\begin{enumerate}[a)]
 \item
 Let $\flag\colon\emptyset=F_0\subsetneq F_1\subsetneq\dots\subsetneq F_r=\groundset$ be a chain and $f\colon\groundset\to\Z_{\geq1}$ be a weight function.
Suppose that $f$ has a set composition type that deconcatenates with respect to $\flag$ (see \cref{ssec:wqsym} for detailed definitions).
Then a basis $B$ of $\matroid$ has $f$-maximum weight if and only if, for every $i$, its intersection with $F_i\setminus F_{i-1}$ is an $f$-maximum-weight basis of the minor $\matroid_i=\matroid|_{F_i}/F_{i-1}$.
Intuitively, the greedy algorithm first fills up $F_1$, then continues inside $F_2$ with $F_1$ already spanned.
In particular, $\matroid$ has a unique maximum-weight basis if and only if every $\matroid_i$ does.
This is a reformulation of results in \cite[Section 2]{ardila_bergman_2005}.
See \Cref{sec:factorization}.

 \item
  A chain $\flag$ of length $r$ consists of flats of a loopless matroid $\matroid$ of rank $r$ exactly when every minor $\matroid_i$ is loopless of rank one. See \Cref{lm:flags}.

 \item
  There are explicit rational numbers $a(\opi)$ such that, for every matroid $\matroid$ on a finite set $\groundset$, the sum of $a(\opi)$ over the $\matroid$-generic set compositions $\opi$ of $\groundset$ equals $1$ if $\rk(\matroid)+|\Loops(\matroid)|=1$, and $0$ otherwise.
 See \Cref{sec:detector}.

\end{enumerate}

\section*{AI disclosure}
The proof was suggested by the AI models Claude Fable~5.1 and GPT 6 Astra Pro, with the AI model GPT-5.6 Sol Pro also used in the process, in response to a description of the problem used in \cite{penaguiao2026chromaticb}.
Gemini 4 Neon was used for literature search.
Claude Code Opus was used for preparing an initial manuscript based,  previous work by the authors as well as the prompt output and  board pictures elaborating the structure.
Claude Code Opus was also used for carrying out edits guided by the authors.
The authors take full responsibility for the correctness, relevance and beauty of the mathematics presented here.

\section*{Acknowledgements}

The AI models used to develop the main proof in this work were accessed via \href{https://math.sciencebench.ai/}{ScienceBench},
we thank Christian Stump for providing the model access through a free research credit grant.
Claude Code was used via a personal subscription from RP.

We would like to thank Yakob Kahane and Alejandro H. Morales for their help in improving the proof of \cref{lm:altsum}.

\section{Preliminaries}\label{sec:prelim}
We briefly recall the necessary notions and set the notation here.
For a more detailed introduction see \cite{penaguiao2026chromaticb} and references therein.

\subsection{Matroids}
We use standard matroid terminology, see, e.g.,  \cite{Oxley2006}.
Throughout, $\groundset$ denotes a nonempty finite set.
A matroid $\matroid$ on $\groundset$ is given by its set of bases $\bases(\matroid)$.
We also use the symmetric exchange property:
\begin{lm}[{\cite[Theorem 2]{brualdi1969comments}}]
 If $B,B'$ are bases of $\matroid$ and $b\in B\setminus B'$, then there is $c\in B'\setminus B$ such that both $(B\setminus\{b\})\cup\{c\}$ and $(B'\setminus\{c\})\cup\{b\}$ are bases.
\end{lm}

We write $\rk_\matroid$ for its rank function, $\rk(\matroid)\coloneqq\rk_\matroid(\groundset)$, and $\cl_\matroid$ for its closure operator,
\begin{equation}
 \cl_\matroid(A)\coloneqq\{e\in \groundset\colon\rk_\matroid(A\cup\{e\})=\rk_\matroid(A)\}\,.
\end{equation}
A \Def{flat} is a set $F$ with $\cl_\matroid(F)=F$.
A \Def{loop} is an element contained in no basis, equivalently an element of $\cl_\matroid(\emptyset)$; we write $\Loops(\matroid)$ for the set of loops, and call $\matroid$ \Def{loopless} if $\Loops(\matroid)=\emptyset$, i.e.\ if $\emptyset$ is a flat.

For $T\subseteq \groundset$, the \Def{restriction} $\matroid|_T$ is the matroid on $T$ with rank function $\rk_\matroid$ restricted to subsets of $T$, and the \Def{contraction} $\matroid/T$ is the matroid on $\groundset\setminus T$ with rank function
\begin{equation}\label{eq:contractionrank}
 \rk_{\matroid/T}(A)=\rk_\matroid(A\cup T)-\rk_\matroid(T)\,,\qquad A\subseteq \groundset\setminus T\,.
\end{equation}
Recall that a subset $J\subseteq \groundset$ is \Def{independent} in a matroid $\matroid$ if and only if $\rk_\matroid(J)=\lvert J\rvert$.
We use the following standard description of the collection of bases of the matroid contraction, see, e.g., \cite[Proposition 3.1.8]{Oxley2006}:
For a matroid $\matroid$ on ground set $\groundset$ and a subset $T\subseteq \groundset$ let $B_T$ be a basis of the restricted matroid $\matroid|_T$. Then
\begin{equation}\label{eq:contractionbases}
 \bases(\matroid/T)=\{ B\subseteq \groundset\setminus T\colon B\cup B_T\in\bases(\matroid)\}\,.
\end{equation}

\subsection{Nested matroids and Hampe's theorem}
A flat $F$ of $\matroid$ is \Def{cyclic} if it is a union of circuits.
A matroid is \Def{nested} if its cyclic flats form a chain under inclusion; this is the definition used by Hampe \cite[Definition~2.8]{hampe_intersection_2017}.
A loopless matroid on a nonempty set has rank at least $1$.
By \cite[Theorem~4.5]{hampe_intersection_2017}, for $1\le r\le\nd$ the number of loopless nested matroids of rank $r$ on $[\nd]$ is the Eulerian number $A_{r-1,\nd}$, the number of permutations of $[\nd]$ with exactly $r-1$ ascents.
Summing over $r$, there are $\nd!$ loopless nested matroids on $[\nd]$ in total.

For $r\ge1$ let $\setofflags_r^\nd$ be the set of chains
\begin{equation*}
 \flag\colon\ \emptyset=F_0\subsetneq F_1\subsetneq\cdots\subsetneq F_r=[\nd]
\end{equation*}
of subsets of $[\nd]$ of length $r$.
For a matroid $\matroid$ of rank $\rk(\matroid)=r$ on ground set $[\nd]$ define the \Def{chain-of-flats vector} $v_\matroid\in\Q^{\setofflags_r^\nd}$ by
\begin{equation*}
 (v_\matroid)_\flag\coloneqq\begin{cases}1, & \text{if } F_1,\dots,F_r \text{ are flats of } \matroid,\\ 0, &\text{otherwise.}\end{cases}
\end{equation*}

\begin{thm}[{\cite[Proposition~3.2]{hampe_intersection_2017}}]\label{thm:hampe}
For every $r\ge1$, the vectors $v_\natroid$, as $\natroid$ ranges over the loopless nested matroids of rank $r$ on $[\nd]$, are linearly independent in $\Q^{\setofflags_r^\nd}$.
\end{thm}

\subsection{Set compositions and word-quasisymmetric functions}\label{ssec:wqsym}
A \Def{set composition} of $\groundset$ is a sequence $\opi=\pi_1|\pi_2|\cdots|\pi_\ell$ of nonempty, pairwise disjoint sets (\Def{blocks}) with union $\groundset$; we write $\opi\models \groundset$ and $\ell(\opi)=\ell$ for the number of blocks.
The \Def{type} of a set composition $\opi$ is defined as $\type(\opi)\coloneqq(\lvert\pi_1\rvert,\dots,\lvert\pi_\ell\rvert)$.
The \Def{last block} of $\opi$ is $\pi_{\ell(\opi)}$.
If $\opi\models \groundset$ and $\otau\models T$ with $\groundset\cap T=\emptyset$, their \Def{concatenation} $\opi|\otau\models \groundset\sqcup T$ lists the blocks of $\opi$ followed by those of $\otau$.
For $\opi=\pi_1|\cdots|\pi_\ell\models \groundset$ and $1\le j\le\ell$, the \Def{upper sets} of $\opi$ are
\begin{equation*}
 U_j(\opi)\coloneqq\pi_j\cup\pi_{j+1}\cup\cdots\cup\pi_\ell\,,
\end{equation*}
so that $\pi_\ell=U_\ell(\opi) \subsetneq U_{\ell-1}(\opi)\subsetneq\cdots\subsetneq  U_1(\opi)=\groundset$ is a  chain.

Let $\flag\colon\emptyset=F_0\subsetneq F_1\subsetneq\cdots\subsetneq F_r=\groundset$ be a chain with $D_i\coloneqq F_i\setminus F_{i-1}$.
We say a set composition $\otau=\tau_1|\tau_2|\dots|\tau_\ell \models\groundset$ \Def{deconcatenates} with respect to the chain $\flag$ if every $F_i$ is an upper set of $\otau$.
Equivalently, every block of $\otau$ lies in a single $D_i$, with the blocks lying in $D_r$ coming first, those lying in $D_{r-1}$ coming just after them, and so on.
Said differently, $\otau$ deconcatenates with respect to $\flag$ if and only if $\otau=\opi_r|\dots|\opi_1$ with $\opi_i\models D_i$ denoting the restriction of $\otau$ to $D_i$.

To $\opi\models \groundset$ we associate the weight function $f_{\opi}\colon \groundset\to\Z_{\ge1}$, $f_{\opi}(e)=j$ for $e\in\pi_j$.
Every $f\colon \groundset\to\Z_{\ge1}$ determines a unique set composition $\opi(f)$, listing the nonempty level sets of $f$ in increasing order of the value; we call it the \Def{set composition type} of $f$.
One can see that the set composition type of $f_{\opi}$ is precisely $\opi$.
It is a known fact, and a consequence of the greedy algorithm (see also \cref{lm:maxbases}), that the maximum-weight bases of $f$ in a matroid depend only on the set composition type of $f$.
A set composition $\opi\models \groundset$ is \Def{$\matroid$-generic} if $f_{\opi}$ is, i.e.\ if $\matroid$ has a unique basis of maximum $f_{\opi}$-weight.

The space $\WQSym_\nd$ of word-quasisymmetric functions of degree $\nd$ has the \Def{monomial basis} $\{\Mnco_{\otau}\}_{\otau\models[\nd]}$, where $\Mnco_{\otau}$ is the sum of all monomials $\bw_{f(1)}\cdots\bw_{f(\nd)}$ with $\opi(f)=\otau$ \cite{bergeron2009hopf}.
Since $\matroid$-genericity of $f$ only depends on $\opi(f)$,
\begin{equation}\label{eq:Psimonomial}
 \ncPsi(\matroid)=\sum_{\substack{\otau\models[\nd]\\ \otau \text{ is }\matroid\text{-generic}}} \Mnco_{\otau}\,.
\end{equation}

\begin{smpl}\label{smpl:U13}
Let $\unifMat1{3}$ be the matroid on $[3]$ whose bases are the singletons.
A set composition $\opi\models[3]$ is $\unifMat1{3}$-generic exactly when its last block is a singleton (the unique element of maximum weight).
Hence the generic ones are $12|3$, $13|2$, $23|1$ and the six compositions into singletons, and $\ncPsi(\unifMat13)=\Mnco_{12|3}+\Mnco_{13|2}+\Mnco_{23|1}+\sum_{\sigma}\Mnco_{\sigma(1)|\sigma(2)|\sigma(3)}$, the last sum being over all permutations $\sigma$ of $[3]$.
\end{smpl}

\section{Maximum-weight bases and factorization along a chain}\label{sec:factorization}

In this section we recall and reprove some basic properties of matroids and their base polytopes, which essentially follow from the greedy algorithm.

\begin{lm}[{\cite[Proposition 1]{ardila_bergman_2005}} reformulated]\label{lm:maxbases}
Let $\matroid$ be a matroid on $\groundset$ and $\opi\models \groundset$ with $\ell$ blocks.
A basis $B$ of $\matroid$ has maximum $f_{\opi}$-weight if and only if
\begin{equation*}
 |B\cap U_j(\opi)|=\rk_\matroid(U_j(\opi))\qquad\text{for all } j=1,\dots,\ell\,.
\end{equation*}
\end{lm}
\begin{proof}
Write $U_j=U_j(\opi)$.
An element of $\pi_j$ lies in exactly $j$ of the sets $U_1,\dots,U_\ell$, so for every basis $B$,
\begin{equation*}
 f_{\opi}(B)=\sum_{j=1}^\ell j\,|B\cap\pi_j|=\sum_{j=1}^\ell|B\cap U_j|\le\sum_{j=1}^\ell\rk_\matroid(U_j)\,,
\end{equation*}
where the inequality holds termwise, since $B\cap U_j$ is independent.
By iteratively enlarging independent sets in $U_j$ for $j=\ell, \ldots , 1$ we construct a basis that attains the upper bound.
Hence the maximum weight is $\sum_j\rk_\matroid(U_j)$, and a basis attains it if and only if all the termwise inequalities are equalities.
\end{proof}

Now fix a matroid $\matroid$ on $\groundset$ and a chain $\flag\colon\emptyset=F_0\subsetneq F_1\subsetneq\cdots\subsetneq F_r=\groundset$ of subsets.
For $1\le i\le r$ let
\begin{equation}\label{eq:defDiMi}
 D_i\coloneqq F_i\setminus F_{i-1}\,,\qquad \matroid_i\coloneqq\matroid|_{F_i}/F_{i-1}\,,
\end{equation}
so $\matroid_i$ is a matroid on $D_i$ and, by \eqref{eq:contractionrank},
\begin{equation}\label{eq:minorrank}
 \rk_{\matroid_i}(A)=\rk_\matroid(F_{i-1}\cup A)-\rk_\matroid(F_{i-1})\qquad\text{for } A\subseteq D_i\,.
\end{equation}

\begin{lm}[{\cite[Proposition 2]{ardila_bergman_2005}} reformulated]\label{lm:factorization}
Let $\otau\models\groundset$ be a set composition that deconcatenates with respect to the chain $\flag$, i.e.,
$\otau =\opi_r|\opi_{r-1}|\cdots|\opi_1$ with  $\opi_i\models D_i$ for $i=1,\dots,r$.

Then the map $B\mapsto(B\cap D_1,\dots,B\cap D_r)$ is a bijection between the bases of $\matroid$ of maximum $f_{\otau}$-weight and the tuples $(B_1,\dots,B_r)$, where each $B_i$ is a basis of $\matroid_i$ of maximum $f_{\opi_i}$-weight.
Its inverse is $(B_1,\dots,B_r)\mapsto B_1\sqcup\cdots\sqcup B_r$.

In particular, $\otau$ is $\matroid$-generic if and only if $\opi_i$ is $\matroid_i$-generic for every $i=1,\dots,r$.
\end{lm}

\begin{rem}
For $r=2$, \cref{lm:factorization} is a consequence of the fact that $\ncPsi$ is a morphism of Hopf monoids \cite[Theorem~B.1]{penaguiao2026chromaticb}, the compatibility with the coproduct coming from \cite[Lemma~B.3]{penaguiao2026chromaticb} applied to matroid base polytopes: the coproduct map $\Delta_{A, B}$ of matroids sends $\matroid$ to $\matroid|_{A}\otimes\matroid/A$, and the coproduct of word-quasisymmetric functions deconcatenates set compositions.
This can iteratively be extended to $r\geq3$.
\cref{lm:factorization} expands on the results from \cite{penaguiao2026chromaticb} by establishing a bijection between the maximum-weight bases of $\matroid$ and the tuples of maximum-weight bases of the minors $\matroid_i$.

Geometrically, \cref{lm:factorization} is the well-known fact that the face of the matroid base polytope of $\matroid$ maximizing a linear functional that is constant on the pieces $D_i$ and decreasing in $i$ is the product of the base polytopes of the minors $\matroid_i$ (see, e.g., \cite{aguiar2023hopf}).

\end{rem}

\begin{proof}[Proof of \cref{lm:factorization}]

Let $B$ be an $f_{\otau}$-maximum-weight basis of $\matroid$ and $B_i\coloneqq B\cap D_i$.

We first check that $B_i$ is a basis in $\matroid_i$.
For $i=1,\dots,r$, $F_i$ is an upper set of $\otau$, so \cref{lm:maxbases} gives $|B\cap F_i|=\rk_\matroid(F_i)$; hence $B\cap F_i$ is a basis of $\matroid|_{F_i}$ for every $i=0,\dots,r$ (trivially for $i=0$).
Since the restriction of $\matroid|_{F_i}$ to $F_{i-1}$ is $\matroid|_{F_{i-1}}$, \eqref{eq:contractionbases} shows that $B_i=(B\cap F_i)\setminus F_{i-1}$ is a basis in $\matroid_i$.

The upper sets of $\otau$ are exactly the sets $F_{i-1}\cup U_k(\opi_i)$, where $1\le i\le r$ and $U_k(\opi_i)$ is a (possibly empty) upper set of $\opi_i$.
From \cref{lm:maxbases} and the assumption that $B$ is an $f_{\otau}$-maximal basis we get
\begin{equation}\label{eq:rkBwithupperset}
 |B\cap(F_{i-1}\cup U_k(\opi_i))|=\rk_\matroid(F_{i-1}\cup U_k(\opi_i))
\end{equation}
for every $i$ and every upper set $U_k(\opi_i)$ of $\opi_i$, including $U_k(\opi_i)=\emptyset$.
In particular,
\begin{equation}\label{eq:rkBwithFi}
 |B\cap F_{i-1}|=\rk_\matroid(F_{i-1})\,.
\end{equation}
Note that, since $U_k(\opi_i)\subseteq D_i=F_i\setminus F_{i-1}$, we have
\begin{align}
 B_i\cap U_k(\opi_i)&=B\cap D_i\cap U_k(\opi_i)=B\cap U_k(\opi_i)\label{eq:Bicapupperset}\\
 F_{i-1}\cap U_k(\opi_i)&=\emptyset\,. \label{eq:Fi-1capupperset}
\end{align}
We  compute
\begin{equation*}
\begin{split}
 \rk_{\matroid_i}(U_k(\opi_i))&\stackrel{\eqref{eq:minorrank}}{=}\rk_{\matroid}(F_{i-1}\cup U_k(\opi_i))-\rk_\matroid(F_{i-1})\\
  &\stackrel{\eqref{eq:rkBwithupperset}}{=}\lvert B\cap(F_{i-1}\cup U_k(\opi_i))\rvert-\rk_\matroid(F_{i-1})\\
  &\stackrel{\eqref{eq:Fi-1capupperset}}{=}\lvert B\cap F_{i-1}\rvert + \lvert B\cap U_k(\opi_i)\rvert-\rk_\matroid(F_{i-1})\\
  &\stackrel{\eqref{eq:rkBwithFi}}{=}\lvert B\cap U_k(\opi_i)\rvert \stackrel{\eqref{eq:Bicapupperset}}{=}\lvert B_i \cap U_k(\opi_i)\rvert\,.
  \end{split}
\end{equation*}
It follows that $B_i$ is $f_{\opi_i}$-maximal by \cref{lm:maxbases}.

Conversely, let each $B_i$ be an $f_{\opi_i}$-maximum basis of $\matroid_i$.
Note that $B_1\in\bases(\matroid|_{F_1})$.
Assuming that $B_1\cup\dots\cup B_{i-1}\in\bases(\matroid|_{F_{i-1}})$ it follows from $B_i\in\bases(\matroid_i)$ and \eqref{eq:contractionbases} that $B_1\cup\dots\cup B_i\in\bases(\matroid|_{F_i})$.
By induction we get that $B=B_1\cup\dots\cup B_r\in\bases(\matroid)$.

To show that the basis $B$ is $f_{\otau}$-maximal we (again) use \Cref{lm:maxbases} and compute
\begin{align*}
 \lvert B\cap (F_{i-1}\cup U_k(\opi_i))\rvert \stackrel{\eqref{eq:Fi-1capupperset}}{=}& \lvert B\cap F_{i-1}\rvert + \lvert B\cap U_k(\opi_i)\rvert\\
 \stackrel{\eqref{eq:Bicapupperset}}{=}& \lvert B_1\cup\dots\cup B_{i-1}\rvert + \lvert B_i\cap U_k(\opi_i)\rvert\\
 \stackrel{}{=}& \rk_\matroid(F_{i-1})+\rk_{\matroid_i}(U_k(\opi_i))\\
 \stackrel{\eqref{eq:minorrank}}{=}& \rk_{\matroid|_{F_i}}(U_k(\opi_i)\cup F_{i-1})
 \stackrel{}{=} \rk_\matroid(U_k(\opi_i)\cup F_{i-1})\,.
\end{align*}

The two maps are clearly inverse to each other.
\end{proof}

\begin{smpl}\label{smpl:factorization}
Let $\matroid$ be the matroid on $[4]$ of rank $2$ whose bases are all two-element subsets except $\{3,4\}$.
\begin{enumerate}[(a)]
 \item\label{it:exflagotau} Let $\flag\colon\emptyset\subsetneq\{2,3\}\subsetneq\{1,2,3,4\}$ and let $\otau=14|2|3$, so $f_{\otau}(1,2,3,4)=(1,2,3,1)$.
  Here $\opi_1=2|3\models\{2,3\}=D_1$ and $\opi_2=14\models \{1,4\}=D_2$, so $\otau$ deconcatenates with respect to $\flag$.
  The minor $\matroid_1=\matroid|_{\{2,3\}}$ has the single basis $\{2,3\}$, and $\matroid_2=\matroid/\{2,3\}$ has rank $0$, with the single basis $\emptyset$.
  Both minors have exactly one basis, so all set compositions are generic, and \cref{lm:factorization} predicts a unique maximum-weight basis $\{2,3\}\sqcup\emptyset=\{2,3\}$.
  Indeed $\{2,3\}$ has weight $5$, while every other basis has weight at most $4$.
  \item Let $\flag$ be as in \ref{it:exflagotau} and $\otau=2|3|14$, then $\otau$ does not deconcatenate with respect to $\flag$ and \cref{lm:factorization} does not apply.
 \item Let $\flag\colon\emptyset\subsetneq\{3,4\}\subsetneq\{1,2,3,4\}$ and let $\otau=12|34$, so $f_{\otau}(1,2,3,4)=(1,1,2,2)$.  Here $D_1=\{3,4\}$ and $D_2=\{1,2\}$ and $\otau$ deconcatenates with respect to $\flag$.
  The minor $\matroid_1=\matroid|_{\{3,4\}}$ has rank $1$ and two bases $\{3\}$, $\{4\}$, and $\matroid_2=\matroid/\{3,4\}$ has rank $1$ and two bases $\{1\}$, $\{2\}$.
  Neither one-block composition is generic, and \cref{lm:factorization} predicts $2\cdot2=4$ maximum-weight bases of $\matroid$, namely the unions $\{1,3\}$, $\{1,4\}$, $\{2,3\}$, $\{2,4\}$.
  Indeed these four bases all have weight $3$, the basis $\{1,2\}$ has weight $2$, and $\{3,4\}$ is not a basis.
  So $\otau$ is not $\matroid$-generic, and the tie in $\matroid$ is exactly the product of the ties in the two minors.
 \item Let $\flag\colon\emptyset\subsetneq\{4\}\subsetneq\{1,2,3,4\}$ and let $\otau=123|4$, so $f_{\otau}(1,2,3,4)=(1,1,1,2)$. Here $D_1=\{4\}$, $D_2=\{1,2,3\}$ and $\otau$ deconcatenates with respect to $\flag$.
  The minor $\matroid_1=\matroid|_{\{4\}}$ has the single basis $\{4\}$, and $\matroid_2=\matroid/\{4\}$ has rank $1$ and bases $\{1\}$, $\{2\}$; the element $3$ is a loop of $\matroid_2$.
  So $\matroid_1$ has one maximum-weight basis and $\matroid_2$ has two, and \cref{lm:factorization} predicts $1\cdot2=2$ maximum-weight bases of $\matroid$, namely $\{1,4\}$ and $\{2,4\}$.
  Indeed these have weight $3$, and all other bases have weight $2$.
  So $\otau$ is not $\matroid$-generic.
\end{enumerate}
\end{smpl}

Lastly, we collect some consequences of the definition of flats along a chain $\flag\colon\emptyset=F_0\subsetneq\cdots\subsetneq F_r=\groundset$.
See also \cite[Theorem 1]{ardila_bergman_2005}.
Recall the notation of $\flag$, $D_i$, $\matroid_i$ as defined in Equation~\eqref{eq:defDiMi} in the beginning of \cref{sec:factorization}.
\begin{lm}\label{lm:flags}
Let $\matroid$ be a matroid on $\groundset$ and $\flag\colon\emptyset=F_0\subsetneq\cdots\subsetneq F_r=\groundset$ a chain.
\begin{enumerate}[(a)]
 \item\label{it:ranksum} $\sum_{i=1}^r\rk(\matroid_i)=\rk(\matroid)$.
 \item\label{it:loops} An element $e\in D_i$ is a loop of $\matroid_i$ if and only if $e\in\cl_\matroid(F_{i-1})$.
 \item\label{it:flagsofflats} Suppose $\matroid$ is loopless.
  Then all $\matroid_i$ are loopless of rank $1$ if and only if $\rk(\matroid)=r$ and $F_1,\dots,F_r$ are flats of $\matroid$.
\end{enumerate}
\end{lm}
\begin{proof}
\ref{it:ranksum} By \eqref{eq:minorrank}, $\rk(\matroid_i)=\rk_\matroid(F_i)-\rk_\matroid(F_{i-1})$, and the sum telescopes to $\rk_\matroid(\groundset)-\rk_\matroid(\emptyset)=\rk(\matroid)$.

\ref{it:loops} By \eqref{eq:minorrank}, $\rk_{\matroid_i}(\{e\})=\rk_\matroid(F_{i-1}\cup\{e\})-\rk_\matroid(F_{i-1})$, which vanishes exactly when $e\in\cl_\matroid(F_{i-1})$.

\ref{it:flagsofflats} Suppose all $\matroid_i$ are loopless of rank $1$.
Then $\rk(\matroid)=r$ by \ref{it:ranksum}.
Let $1\le i\le r$ and $e\notin F_i$, say $e\in D_j$ with $j>i$.
Since $\matroid_j$ is loopless, \ref{it:loops} gives $e\notin\cl_\matroid(F_{j-1})\supseteq\cl_\matroid(F_i)$.
Hence $\cl_\matroid(F_i)=F_i$.

Conversely, suppose $\rk(\matroid)=r$ and all $F_i$ are flats.
Since $\matroid$ is loopless, $F_0=\emptyset$ is a flat as well.
For $1\le i\le r$, an element of $D_i$ lies outside the flat $F_{i-1}$, so $\rk_\matroid(F_i)>\rk_\matroid(F_{i-1})$.
The $r$ strict increases $0=\rk_\matroid(F_0)<\cdots<\rk_\matroid(F_r)=r$ force $\rk_\matroid(F_i)-\rk_\matroid(F_{i-1})=1$, i.e.\ $\rk(\matroid_i)=1$.
Finally, $\cl_\matroid(F_{i-1})\cap D_i=F_{i-1}\cap D_i=\emptyset$, so $\matroid_i$ is loopless by \ref{it:loops}.
\end{proof}

\section{A local detector for rank one}\label{sec:detector}

In this section we consider a poset defined for a basis in a matroid.
We use this poset to prove identities of alternating sums over this poset (\Cref{lm:altsum}) and finally define a map (Equation~\eqref{eq:defdelta}) that detects matroids of rank one without loops and matroids consisting of only one loop (\Cref{prop:detector}).

For a matroid $\matroid$ on $\groundset$ and a basis $B$, an \Def{exchange pair} is a pair $(b,c)$ with $b\in B$, $c\notin B$ and $(B\setminus\{b\})\cup\{c\}\in\bases(\matroid)$.
We encode the exchange pairs of $B$ in a \Def{basis exchange poset} $\poset_B$ on $\groundset$: declare $c<_{\poset_B}b$ whenever $(b,c)$ is an exchange pair.
This poset has been used in \cite{billera2009quasisymmetric} and \cite{luoto_matroid-friendly_2008}.
Since every such relation has its larger element in $B$ and its smaller element outside $B$, no two relations can be chained, so this is already a partial order (of height at most one).
The Hasse diagram of $\poset_B$ is the \Def{dependence graph} of $\matroid$ with respect to $B$, studied by Krogdahl \cite{krogdahl1977dependence}: the bipartite graph on $B\sqcup(\groundset\setminus B)$ joining $b\in B$ and $c\notin B$ whenever $b$ lies in the fundamental circuit of $c$ with respect to $B$, equivalently whenever $(b,c)$ is an exchange pair.
We start by recording and proving two elementary observations.
\begin{lm}\label{lm:localopt}
Let $\matroid$ be a matroid on $\groundset$, $B$ a basis and $f\colon \groundset\to\R$.
Then $B$ is the unique basis of maximum $f$-weight if and only if $f(b)>f(c)$ for every exchange pair $(b,c)$ of $B$.
\end{lm}
\begin{proof}
If $f(b)\le f(c)$ for some exchange pair, the basis $(B\setminus\{b\})\cup\{c\}\neq B$ has weight at least $f(B)$.
Conversely, suppose $f(b)>f(c)$ for all exchange pairs, and assume for a contradiction that there is a basis $B'\neq B$ with $f(B')\ge f(B)$.
Choose such a $B'$ with $|B\setminus B'|$ as small as possible.
Pick $b\in B\setminus B'$; by the symmetric exchange property there is $c\in B'\setminus B$ such that $(B\setminus\{b\})\cup\{c\}$ and $B''\coloneqq(B'\setminus\{c\})\cup\{b\}$ are bases.
Then $(b,c)$ is an exchange pair of $B$, so $f(b)>f(c)$ and
\begin{equation*}
 f(B'')=f(B')+f(b)-f(c)>f(B')\ge f(B)\,.
\end{equation*}
In particular $B''\neq B$.
But $|B\setminus B''|=|B\setminus B'|-1$, contradicting the choice of $B'$.
\end{proof}
\begin{lm}\label{lm:maxelements}
The maximal elements of $\poset_B$ are exactly the elements of $B\sqcup\Loops(\matroid)$.
\end{lm}
\begin{proof}
Elements of $B$ are never smaller than anything, so they are maximal.
A loop $c$ lies in no basis, so it is in no exchange pair, and is isolated in $\poset_B$, hence maximal.
Finally, let $c\notin B$ be a non-loop.
Then $B\cup\{c\}$ contains a unique circuit $C$, and $c\in C$ \cite{Oxley2006}.
Since $c$ is not a loop, $C\neq\{c\}$, so there is some $b\in C\setminus\{c\}\subseteq B$.
The set $(B\setminus\{b\})\cup\{c\}$ is contained in $B\cup\{c\}$ and does not contain the only circuit $C$ of $B\cup\{c\}$, so it is independent; having $|B|$ elements, it is a basis.
Hence $c<_{\poset_B}b$, and $c$ is not maximal.
\end{proof}

Say that $\opi\models \groundset$ is \Def{strictly compatible} with a poset $\poset$ on $\groundset$ if $x<_\poset y$ implies $f_{\opi}(x)<f_{\opi}(y)$, and we write $\opi\models \poset$.
Equivalently, $\opi$ is strictly compatible with a poset $\poset$ if and only if $f_{\opi}\colon\groundset\to\Z_{\geq 1}$ is strictly $\poset$-order preserving.
By \cref{lm:localopt}, for a basis $B$ of $\matroid$,
\begin{equation}\label{eq:genericposet}
 \opi \text{ is } \matroid\text{-generic with maximum-weight basis } B \iff \opi \models \poset_B\,.
\end{equation}

\Cref{lm:altsum} below can be proven in an elementary way without the use of generating functions;
in fact, this is what the AI agents suggested.
We instead give a proof via order polynomials, in a slightly more general setting.
Similar functions appear in \cite{liu_p-partition_2020,liu_p-partitions_2021}.

\begin{lm}\label{lm:altsum}
Let $\poset$ be a poset on a set $\groundset$ with $|\groundset|\ge1$, and let $\Max(\poset)$ be its set of maximal elements.
Then
\begin{equation*}
 \sum_{\opi\models \poset} (-1)^{|\groundset|-\ell(\opi)}\,|\pi_{\ell(\opi)}| \;=\;\begin{cases}1,&\text{if }|\Max(\poset)|=1,\\0,&\text{otherwise.}\end{cases}
\end{equation*}
\end{lm}
\begin{proof}
Consider the following polynomial:
\begin{equation}
 \poly_\poset(t,n) = \sum_{\opi\models \poset} \binom{n}{\ell(\opi)} t^{\lvert\pi_{\ell(\opi)}\rvert}\,.
\end{equation}

We specialize $t=1$ and we get
\begin{equation}
 \poly_\poset(1, n) = \sum_{\opi\models \poset} \binom{n}{\ell(\opi)} =\Omega^\circ (\poset,n)\,,
\end{equation}
where $\Omega^\circ (\poset,n)$ denotes the strict order polynomial.
From combinatorial reciprocity (see, e.g., \cite[Corollary~3.15.12]{stanley2011enumerative}) and the fact that there is exactly one constant order-preserving map we get
\begin{equation}\label{eq:step1}
\begin{split}
 \sum_{\opi \models \poset}(-1)^{\ell(\opi)} =  \sum_{\opi\models \poset} \binom{-1}{\ell(\opi)}= \poly_\poset(1,-1) =\Omega^\circ (\poset,-1)
 = (-1)^{|\groundset|}\Omega(\poset,1)=(-1)^{|\groundset|}\,.
 \end{split}
\end{equation}

Let $\opi$ be a set composition strictly compatible with $\poset$ and write $\opi=\opi'|A$ with last block $A$. Then $A$ is a non-empty subset of the maximal elements $\Max(\poset)$ of $\poset$, $\opi'$ is a set composition strictly compatible with $\poset|_{\groundset\setminus A}$ and $\ell(\opi)=\ell(\opi')+1$.
We now specialize $n=-1$ in $ \poly_\poset(t,n)$, group the sum according to $A$ and compute:
\begin{equation}\label{eq:polynom}
\begin{split}
 \poly_\poset(t, -1) &= \sum_{\opi\models \poset} (-1)^{\ell(\opi)} t^{\lvert\pi_{\ell(\opi)}\rvert}\\
  &= \sum_{\emptyset\neq A\subseteq \Max(\poset)} (-1) t^{\lvert A\rvert} \sum_{\opi'\models \poset|_{\groundset\setminus A}}(-1)^{\ell(\opi')}\\
  &\overset{\eqref{eq:step1}}= \sum_{\emptyset\neq A\subseteq \Max(\poset)} (-1) t^{\lvert A\rvert} (-1)^{\lvert \groundset\rvert-\lvert A\rvert}\\
  &=(-1)^{\lvert \groundset\rvert-1}\sum_{\emptyset\neq A\subseteq \Max(\poset)}  (-t)^{\lvert A\rvert} \\
  &=(-1)^{\lvert \groundset\rvert-1}((1-t)^{\lvert\Max(\poset)\rvert}-1)\,.
  \end{split}
\end{equation}
Taking the derivative with respect to $t$ on both ends of \eqref{eq:polynom} gives:
\begin{equation}
 \sum_{\opi\models \poset} (-1)^{\ell(\opi)} \,\lvert\pi_{\ell(\opi)}\rvert\, t^{\lvert\pi_{\ell(\opi)}\rvert-1}=\frac{\partial \poly_\poset(t,-1)}{\partial t}= (-1)^{\lvert \groundset\rvert}\,\lvert\Max(\poset)\rvert\,(1-t)^{\lvert\Max(\poset)\rvert-1}\,.
\end{equation}
Evaluating $t=1$ and multiplying by $(-1)^{\lvert\groundset\rvert}$ yields the result:
\begin{equation}
 \sum_{\opi\models \poset} (-1)^{\lvert \groundset\rvert-\ell(\opi)} \lvert\pi_{\ell(\opi)}\rvert = \lvert\Max(\poset)\rvert\cdot0^{\lvert\Max(\poset)\rvert-1}=\begin{cases}1,&\text{if }|\Max(\poset)|=1,\\0,&\text{otherwise.}                                                                                                                          \end{cases}
\end{equation}
\end{proof}

For a nonempty finite set $\groundset$  and $\opi\models \groundset$ define
\begin{equation}\label{eq:defa}
 a(\opi)\coloneqq\frac{(-1)^{|\groundset|-\ell(\opi)}\,|\pi_{\ell(\opi)}|}{|\groundset|}\in\Q\,,
\end{equation}
and for a matroid $\matroid$ on $\groundset$ define the \Def{rank-one detector}
\begin{equation}\label{eq:defdelta}
 \delta_\groundset(\matroid)\coloneqq\sum_{\substack{\opi\models \groundset\\ \opi\text{ is }\matroid\text{-generic}}}a(\opi)\,.
\end{equation}

\begin{prop}\label{prop:detector}
For every matroid $\matroid$ on a nonempty finite set $\groundset$,
\begin{equation*}
 \delta_\groundset(\matroid)=\begin{cases}1,&\text{if }\rk(\matroid)+|\Loops(\matroid)|=1,\\0,&\text{otherwise.}\end{cases}
\end{equation*}
In particular, if $\matroid$ is loopless, then $\delta_\groundset(\matroid)=1$ if $\rk(\matroid)=1$ and $\delta_\groundset(\matroid)=0$ otherwise.
\end{prop}
\begin{proof}
Every generic $\opi$ has exactly one maximum-weight basis, so by \eqref{eq:genericposet} the generic compositions are partitioned according to this basis:
\begin{equation*}
 \delta_\groundset(\matroid)=\frac{1}{\lvert\groundset\rvert}\sum_{B\in\bases(\matroid)}\ \sum_{\opi \models \poset_B}(-1)^{|\groundset|-\ell(\opi)}|\pi_{\ell(\opi)}|\,.
\end{equation*}
By \cref{lm:altsum,lm:maxelements}, the inner sum is $1$ if $|B|+|\Loops(\matroid)|=\rk(\matroid)+|\Loops(\matroid)|=1$ and $0$ otherwise, independently of $B$.
Hence
\begin{equation*}
 \delta_\groundset(\matroid)=\begin{cases}
                              \frac{|\bases(\matroid)|}{|\groundset|}\,, & \text{ if }\rk(\matroid)+|\Loops(\matroid)|=1\\
                              0\,,&\text{ otherwise.}
                             \end{cases}
\end{equation*}
Assume $\rk(\matroid)+|\Loops(\matroid)|=1$.
There are two cases: either $\matroid$ is loopless of rank $1$, and then every element of $\groundset$ forms a basis, so $|\bases(\matroid)|=|\groundset|$; or $\matroid$ has rank $0$ and a single loop, and then $\groundset$ consists of that loop, so $|\groundset|=1=|\bases(\matroid)|$.
In both cases $\delta_\groundset(\matroid)=1$.
\end{proof}

\begin{smpl}\label{smpl:delta}
For $\groundset=[3]$: $a(123)=1$; $a(A|B)=-|B|/3$ for the six compositions with two blocks; and $a(\opi)=1/3$ for the six compositions into singletons.
For $\matroid=\unifMat13$, \cref{smpl:U13} gives $\delta_{[3]}(\matroid)=3\cdot(-\tfrac13)+6\cdot\tfrac13=1$.

For $\matroid=\unifMat23$, the generic compositions are $1|23$, $2|13$, $3|12$ (with unique maximum basis the last block) and the six compositions into singletons, so $\delta_{[3]}(\matroid)=3\cdot(-\tfrac23)+6\cdot\tfrac13=0$, as predicted by \cref{prop:detector}.

For $\groundset=[4]$ there are $75$ set compositions.
The value $a(\opi)$ only depends on the type $\type(\opi)=(|\pi_1|,\dots,|\pi_\ell|)$, and so does genericity for the uniform matroids: $\opi$ is $\unifMat k4$-generic if and only if the $k$ largest values of $f_{\opi}$ are separated from the others, that is, if and only if $k=|U_j(\opi)|$ for some $j$.
\cref{table:exampledeltas} lists, for each type, the number of set compositions of that type, the value of $a$, and how many of them are generic for $\unifMat14$, $\unifMat24$, $\unifMat34$, and for the matroid $\mathrm L$ on $[4]$ with bases $\{1\},\{2\},\{3\}$, in which $4$ is a loop.
\begin{table}
\begin{center}
\renewcommand{\arraystretch}{1.2}
\begin{tabular}{c|c|c|cccc}
 type & \#$\opi$ & $a(\opi)$ & $\unifMat14$ & $\unifMat24$ & $\unifMat34$ & $\mathrm L$\\
 \hline
 $(4)$ & $1$ & $-1$ & $0$ & $0$ & $0$ & $0$\\
 $(1,3)$ & $4$ & $\tfrac34$ & $0$ & $0$ & $4$ & $0$\\
 $(2,2)$ & $6$ & $\tfrac12$ & $0$ & $6$ & $0$ & $3$\\
 $(3,1)$ & $4$ & $\tfrac14$ & $4$ & $0$ & $0$ & $3$\\
 $(1,1,2)$ & $12$ & $-\tfrac12$ & $0$ & $12$ & $12$ & $6$\\
 $(1,2,1)$ & $12$ & $-\tfrac14$ & $12$ & $0$ & $12$ & $9$\\
 $(2,1,1)$ & $12$ & $-\tfrac14$ & $12$ & $12$ & $0$ & $12$\\
 $(1,1,1,1)$ & $24$ & $\tfrac14$ & $24$ & $24$ & $24$ & $24$\\
 \hline
 \multicolumn{3}{r|}{$\delta_{[4]}$} & $1$ & $0$ & $0$ & $0$
\end{tabular}
\end{center}
 \caption{The table lists all the types of set compositions on a ground set of size four, how many set compositions of this type exist, their value $a$ as defined in Equation~\eqref{eq:defa}, and how many set compositions of each type are generic for the matroids in \cref{smpl:delta}.
 The last row computes $\delta(\matroid)$ by summing over the rows and for each row multiplying the entry in the column computing the value $a$ with the number of generic set compositions of that type (the last group of columns).
 }\label{table:exampledeltas}
\end{table}

For instance, for $\unifMat24$ the generic compositions are those of types $(2,2)$, $(1,1,2)$, $(2,1,1)$ and $(1,1,1,1)$, and
\begin{equation*}
 \delta_{[4]}(\unifMat24)=6\cdot\tfrac12+12\cdot\big(-\tfrac12\big)+12\cdot\big(-\tfrac14\big)+24\cdot\tfrac14=3-6-3+6=0\,,
\end{equation*}
while for $\unifMat14$ we get $4\cdot\tfrac14+12\cdot(-\tfrac14)+12\cdot(-\tfrac14)+24\cdot\tfrac14=1-3-3+6=1$.
The matroid $\mathrm L$ has rank $1$ but also a loop, so $\rk(\mathrm L)+|\Loops(\mathrm L)|=2$ and $\delta_{[4]}(\mathrm L)=\tfrac32+\tfrac34-3-\tfrac94-3+6=0$: the loop is detected, although $\mathrm L$ has rank one.
In agreement with \cref{prop:detector}, only $\unifMat14$, the unique loopless matroid of rank one on $[4]$, gives a nonzero value.
\end{smpl}

\section{Chain functionals}\label{sec:chain}

Recall that $\otau$ deconcatenates with respect to $\flag$ if and only if $\otau=\opi_r|\dots|\opi_1$ with $\opi_i\models D_i$ denoting the restriction of $\otau$ to $F_i\setminus F_{i-1}$ (note the reversed order).
We can now define the linear functional $\Lambda_\flag\colon\WQSym_\nd\to\Q$ on the monomial basis by
\begin{equation}\label{eqdefin:chainfunctional}
 \Lambda_\flag(\Mnco_{\otau})\coloneqq \begin{cases}
 \displaystyle\prod_{i=1}^r a(\opi_i)\,, &\text{if } \otau=\opi_r|\dots|\opi_1 \text{ deconcatenates with respect to } \flag \\
 0\,, & \text{otherwise,}\end{cases}
\end{equation}
for $\otau\models[\nd]$, and extend linearly.

\begin{smpl}
 For example, for the chain  $\flag\colon\{2,3\}\subsetneq[4]$ we have
 \begin{align*}
  \Lambda_\flag(\Mnco_{14|23})&=a(23)\,a(14)=(-1)\cdot(-1)=1\\
  \Lambda_\flag(\Mnco_{1|4|2|3})&=a(2|3)\,a(1|4)=\tfrac12\cdot\tfrac12=\tfrac14\,,
 \end{align*}
 while $\Lambda_\flag(\Mnco_{134|2})=0$.
Note the order of concatenation: the blocks coming from $D_r$ are placed \emph{first}.
In this way, $\Lambda_\flag(\Mnco_{23|14})=0$.
This order matches \cref{lm:factorization}, in which the greedy algorithm for maximum-weight bases treats $F_1$ first and agrees with the one given by the coproduct in $\WQSym$.
\end{smpl}

Before proving the crucial \cref{prop:chainfunctional} we give a big picture explanation of how the linear map $\Lambda_\flag$ interacts with the coproduct and the rank-one detector from \cref{prop:detector}.

The functional $\Lambda_\flag$ fits in a commutative diagram, which organizes the proof of \cref{prop:chainfunctional} below.
For a finite set $\groundset$, write $\Mat_\groundset$ for the vector space spanned by the matroids on $\groundset$, and $\WQSym_\groundset$ for the vector space with basis $\{\Mnco_{\opi}\}_{\opi\models \groundset}$;%
\footnote{Here we are using this only for sets $\groundset\subseteq[\nd]$ and hence we naturally have a linear order on $\groundset=\{i_1<\dots< i_n\}$ to define monomials $\bw_{f(i_1)}\bw_{f(i_2)}\dots\bw_{f(i_n)}$ for maps $f\colon \groundset\to\Z_{\geq1}$.}
the map $\ncPsi\colon\Mat_\groundset\to\WQSym_\groundset$ is defined by \eqref{eq:Psimonomial}, verbatim.

For a chain $\flag\colon\emptyset=F_0\subsetneq F_1\subsetneq\cdots\subsetneq F_r=[\nd]$ recall that $D_i\coloneqq F_i\setminus F_{i-1}$.
Consider the linear maps
\begin{align*}
 \Delta_{D_r,\dots,D_1}\colon\Mat_{[\nd]}&\to\bigotimes_{i=r}^1\Mat_{D_i}\,, & \matroid&\mapsto\matroid_r\otimes\cdots\otimes\matroid_1\,,\quad \matroid_i=\matroid|_{F_i}/F_{i-1}\,,\\
 \Delta_{D_r,\dots,D_1}\colon\WQSym_{[\nd]}&\to\bigotimes_{i=r}^1\WQSym_{D_i}\,, & \Mnco_{\otau}&\mapsto\begin{cases}\Mnco_{\opi_r}\otimes\cdots\otimes\Mnco_{\opi_1}\,, & \begin{array}[t]{@{}l@{}}\text{if } \otau=\opi_r|\cdots|\opi_1\\ \text{with } \opi_i\models D_i\,,\end{array}\\[3ex] 0\,, & \text{otherwise,}\end{cases}\\
 a\colon\WQSym_\groundset&\to\Q\,, & \Mnco_{\opi}&\mapsto a(\opi)\,.
\end{align*}
The two maps $\Delta_{D_r,\dots,D_1}$ are iterated coproducts, of the Hopf monoid of matroids (restriction and contraction) and of word-quasisymmetric functions (deconcatenation), respectively\footnote{We use the same notation $\Delta_{D_r,\dots,D_1}$ for the coproduct on matroids and word-quasisymmetric functions. This is justified since the two operations are related by the Hopf algebra morphism $\ncPsi$ and it will be clear from the context which Hopf algebra we are operating on.}.
Identifying $\Q^{\otimes r}$ with $\Q$ via multiplication, we can write Equation~\eqref{eqdefin:chainfunctional} as $\Lambda_\flag=a^{\otimes r}\circ\Delta_{D_r,\dots,D_1}\colon \WQSym\to\Q$,
and by definition $\delta_\groundset=a\circ\ncPsi$ on matroids on $\groundset$.
We obtain the diagram
\begin{equation}\label{eq:diagram}
\begin{tikzcd}[column sep=huge]
 & & \Q \\
 \WQSym_{[\nd]} \arrow[r, "\Delta_{D_r,\dots,D_1}"'] \arrow[urr, "\Lambda_\flag"] & \displaystyle\bigotimes_{i=1}^r\WQSym_{D_i} \arrow[ur, "a^{\otimes r}"' near start] & \\
 \Mat_{[\nd]} \arrow[u, "\ncPsi"] \arrow[r, "\Delta_{D_r,\dots,D_1}"'] & \displaystyle\bigotimes_{i=1}^r\Mat_{D_i} \arrow[u, "\ncPsi^{\otimes r}"'] \arrow[uur, "\bigotimes_i\delta_{D_i}"'] &
\end{tikzcd}
\end{equation}
The top and right triangles commute by the definitions just given, and the square commutes because $\ncPsi$ is a morphism of Hopf monoids \cite[Theorem~B.1]{penaguiao2026chromaticb}, and in particular preserves the coproduct (see \cite[Lemma~B.3]{penaguiao2026chromaticb}, applied to matroid base polytopes).
It also follows directly from \cref{lm:factorization}: for a matroid $\matroid$ on $[\nd]$,
\begin{align*}
 \Delta_{D_r,\dots,D_1}(\ncPsi(\matroid))
 =\sum_{\substack{\opi_1,\dots,\opi_r\\ \opi_r|\cdots|\opi_1\ \matroid\text{-generic}}}\Mnco_{\opi_r}\otimes\cdots\otimes\Mnco_{\opi_1}
 &=\bigotimes_{i=r}^1\sum_{\substack{\opi_i\models D_i\\ \matroid_i\text{-generic}}}\Mnco_{\opi_i}\\
 &=\ncPsi(\matroid_r)\otimes\cdots\otimes\ncPsi(\matroid_1)\,.
\end{align*}
Following $\matroid\in\Mat_{[\nd]}$ along the two paths to $\Q$ gives exactly Equation~\eqref{eq:product} below.

\begin{prop}\label{prop:chainfunctional}
Let $\matroid$ be a matroid on $[\nd]$ and $\flag\in\setofflags_r^\nd$.
Let $D_i$ and $\matroid_i$ be as in \cref{sec:factorization}.

Then
\begin{equation}\label{eq:product}
 \Lambda_\flag(\ncPsi(\matroid))=\prod_{i=1}^r\delta_{D_i}(\matroid_i)\,.
\end{equation}
If moreover $\matroid$ is loopless, then:
\begin{enumerate}[(a)]
 \item\label{it:higherrank} if $\rk(\matroid)>r$, then $\Lambda_\flag(\ncPsi(\matroid))=0$;
 \item\label{it:samerank} if $\rk(\matroid)=r$, then $\Lambda_\flag(\ncPsi(\matroid))=(v_\matroid)_\flag$.
\end{enumerate}
\end{prop}
\begin{proof}
By the linearity of $\Lambda_\flag$ and \eqref{eq:Psimonomial},
\begin{equation}
 \Lambda_\flag(\ncPsi(\matroid))=\sum_{\substack{\otau\\\matroid\text{-generic}}}\Lambda_\flag(\Mnco_{\otau})\,.
\end{equation}
By \eqref{eqdefin:chainfunctional}, only the $\otau$ of the form $\opi_r|\cdots|\opi_1$ with $\opi_i\models D_i$ contribute, so by \cref{lm:factorization}
\begin{align*}
 \Lambda_\flag(\ncPsi(\matroid))
 =\sum_{\substack{\opi_1\models D_1,\dots,\opi_r\models D_r\\ \opi_r|\cdots|\opi_1\ \matroid\text{-generic}}}\ \prod_{i=1}^ra(\opi_i)
 &=\sum_{\substack{\opi_1\models D_1,\dots,\opi_r\models D_r\\ \opi_i\ \matroid_i\text{-generic for all } i}}\ \prod_{i=1}^ra(\opi_i)\\
 &=\prod_{i=1}^r\sum_{\substack{\opi_i\models D_i\\ \opi_i\text{ is }\matroid_i\text{-generic}}}a(\opi_i)
 \ =\ \prod_{i=1}^r \delta_{D_i}(\matroid_i)
\end{align*}
where the third equality is distributivity and the last one is the definition of $\delta_{D_i}$.

Now let $\matroid$ be loopless and suppose the product \eqref{eq:product} is nonzero.
Then $\delta_{D_i}(\matroid_i)=1$ for every $i$, so by \cref{prop:detector} each $\rk(\matroid_i)\le1$, and \cref{lm:flags}\ref{it:ranksum} gives $\rk(\matroid)\le r$.
This proves \ref{it:higherrank}.

For \ref{it:samerank}, let $\rk(\matroid)=r$.
If the product is nonzero, then as before $\rk(\matroid_i)\le 1$ for all $i$, and since these ranks sum to $r$, all $\matroid_i$ have rank $1$; \cref{prop:detector} then forces them to be loopless, so by \cref{lm:flags}\ref{it:flagsofflats} all $F_i$ are flats.
Conversely, if all $F_i$ are flats, \cref{lm:flags}\ref{it:flagsofflats} and \cref{prop:detector} show that every factor equals $1$.
\end{proof}

\begin{smpl}
Let $\matroid=\unifMat23$ on $[3]$, of rank $2$.
Its flats are $\emptyset$, the singletons and $[3]$.
For the chain $\{1\}\subsetneq[3]$ we get $\matroid_1=\matroid|_{\{1\}}$, loopless of rank one, and $\matroid_2=\matroid/\{1\}$, the matroid on $\{2,3\}$ whose bases are $\{2\}$ and $\{3\}$; so
\begin{equation*}
 \Lambda_{\flag}(\ncPsi(\matroid)) = a(1)\cdot a(2|3)+a(1)\cdot a(3|2)=1\cdot \frac 1 2 + 1\cdot  \frac 1 2 =1 = 1\cdot 1 =\delta_{\{1\}}(\matroid_1)\cdot\delta_{\{2,3\}}(\matroid_2)\,.
\end{equation*}
For the chain $\{1,2\}\subsetneq[3]$ we get $\matroid_1=\matroid|_{\{1,2\}}$ of rank $2$, so $\delta_{D_1}(\matroid_1)=0$ and $\Lambda_{\flag}(\ncPsi(\unifMat{2}{3}))=0$, reflecting that $\{1,2\}$ is not a flat.
\end{smpl}

\begin{smpl}\label{smpl:unif2d}\label{smpl:lowrank}
More generally, let $\matroid=\unifMat2\nd$ on $[\nd]$, and let $\flag$ be a chain with pieces of sizes $c=(|D_1|,\dots,|D_r|)$.
Since $\matroid$ is uniform, the minors only depend on these sizes:
\begin{itemize}
 \item $\matroid_1=\matroid|_{D_1}$ is $\unifMat1{D_1}$ if $|D_1|=1$ and $\unifMat2{D_1}$ otherwise, so $\delta_{D_1}(\matroid_1)=1$ if and only if $|D_1|=1$;
 \item if $|D_1|=1$, then $\matroid_2=\matroid|_{F_2}/F_1=\unifMat1{D_2}$ is loopless of rank one, so $\delta_{D_2}(\matroid_2)=1$;
 \item for $i\ge3$ we have $|F_{i-1}|\ge2$, so $F_{i-1}$ spans $\matroid$ and every element of $D_i$ is a loop of $\matroid_i$; hence $\delta_{D_i}(\matroid_i)=1$ if and only if $|D_i|=1$.
\end{itemize}
By \eqref{eq:product}, $\Lambda_{\flag}(\ncPsi(\unifMat2\nd))=1$ if $|D_1|=1$ and $|D_i|=1$ for all $i\ge3$, and $\Lambda_{\flag}(\ncPsi(\unifMat2\nd))=0$ otherwise.
The tables in \cref{table:Unif2d} list every size vector $c$ for $\nd=4$ and $\nd=5$,
and the value $\Lambda_{\flag}(\ncPsi(\unifMat2\nd))$.
\begin{table}
 \begin{center}
\renewcommand{\arraystretch}{1.15}
\begin{tabular}{c|c}
 \multicolumn{2}{c}{$\unifMat24$}\\
 $c$ & $\Lambda_{\flag}$\\
 \hline
 $(4)$  & $0$\\
 $(1,3)$ & $1$\\
 $(2,2)$  & $0$\\
 $(3,1)$ & $0$\\
 $(1,1,2)$  & $0$\\
 $(1,2,1)$  & $1$\\
 $(2,1,1)$  & $0$\\
 $(1,1,1,1)$  & $1$
\end{tabular}
\qquad
\begin{tabular}{c|c}
 \multicolumn{2}{c}{$\unifMat25$}\\
 $c$ & $\Lambda_{\flag}$\\
 \hline
 $(5)$ &  $0$\\
 $(1,4)$  & $1$\\
 $(2,3)$  & $0$\\
 $(3,2)$  & $0$\\
 $(4,1)$ & $0$\\
 $(1,1,3)$  & $0$\\
 $(1,2,2)$  & $0$\\
 $(1,3,1)$  & $1$
\end{tabular}
\qquad
\begin{tabular}{c|c}
 \multicolumn{2}{c}{$\unifMat25$ (cont.)}\\
 $c$ & $\Lambda_{\flag}$\\
 \hline
 $(2,1,2)$  & $0$\\
 $(2,2,1)$  & $0$\\
 $(3,1,1)$  & $0$\\
 $(1,1,1,2)$  & $0$\\
 $(1,1,2,1)$  & $0$\\
 $(1,2,1,1)$ & $1$\\
 $(2,1,1,1)$  & $0$\\
 $(1,1,1,1,1)$  & $1$
\end{tabular}
\end{center}
\caption{Computations for \cref{smpl:unif2d}:  every size vector $c$ of chains of sets of size $\nd=4$ and $\nd=5$,
and the value $\Lambda_{\flag}(\ncPsi(\unifMat2\nd))$}\label{table:Unif2d}
\end{table}
The two regimes of \cref{prop:chainfunctional} are visible:
\begin{itemize}
 \item for $r=1<\rk(\unifMat{2}{5})=2$ the value is $0$, as in \cref{prop:chainfunctional}\ref{it:higherrank};
 \item for $r=2=\rk(\unifMat{2}{5})$ the value is $1$ exactly for the $\nd$ chains $\{j\}\subsetneq[\nd]$, whose first set is a flat of rank one, and $0$ for the chains whose first set has size at least $2$, which is not a flat; this is the vector $v_{\unifMat2\nd}$, as in \cref{prop:chainfunctional}\ref{it:samerank}.
\end{itemize}
Nothing is claimed in \cref{prop:chainfunctional} when $\rk(\matroid)<r$.
Note that the lattice of flats of a matroid $\matroid$ is ranked by its rank function $\rk_\matroid$.
Hence a chain $\flag$ of length $r>\rk(\matroid)$ can never be a chain of flats of $\matroid$.
Nevertheless, we observe both zero and nonzero values of  $\Lambda_\flag(\ncPsi(\unifMat{2}{4}))$ and $\Lambda_\flag(\ncPsi(\unifMat{2}{5}))$ for chains $\flag$ of length $r>2=\rk(\unifMat{2}{4})=\rk(\unifMat{2}{5})$.
\end{smpl}

\section{Proof of the main theorem}\label{sec:main}

\begin{proof}[Proof of \cref{thm:main}]
Suppose that
\begin{equation*}
 \sum_{\natroid \text{ loopless nested matroid on }[\nd]}\lambda_\natroid\ncPsi(\natroid)=0\,,
\end{equation*}
with not all $\lambda_\natroid$ zero.
Let $r$ be the smallest rank of a loopless nested matroid $\natroid$ with $\lambda_\natroid\neq0$.
Apply $\Lambda_\flag$ for $\flag\in\setofflags_r^\nd$.
Summands corresponding to loopless nested matroids of rank greater than $r$ vanish by \cref{prop:chainfunctional}\ref{it:higherrank}, summands corresponding to loopless nested matroids of rank less than $r$ have $\lambda_\natroid=0$ by assumption, and the summands corresponding to loopless nested matroids of rank equal to $r$ give $(v_\natroid)_{\flag}$ by \cref{prop:chainfunctional}\ref{it:samerank}.
Hence
\begin{equation*}
\sum_{\substack{\natroid \text{ loopless nested matroid on }\\ \rk(\natroid)=r}}\lambda_\natroid\,v_\natroid=0\qquad\text{in }\Q^{\setofflags_r^\nd}\,,
\end{equation*}
and \cref{thm:hampe} gives $\lambda_\natroid=0$ for all $\natroid$ of rank $r$, contradicting the choice of $r$.

The dimension result follows from the fact that $\nd!$ loopless nested matroids span $(\im\ncPsi)_\nd$, see \cite[Theorem~1.1]{penaguiao2026chromaticb}.
\end{proof}

\begin{rem}
Hampe \cite{hampe_intersection_2017} showed that \cref{thm:hampe} holds over any field, so \cref{thm:main} holds for any field of characteristic $0$ or $p>\nd$.
\end{rem}

\begin{cor}
The kernel of $\ncPsi\colon\Mat_\nd\to\WQSym_\nd$ is spanned by the valuative relations and the loop-coloop relations as described in \cite{penaguiao2026chromaticb}.
\end{cor}

\printbibliography

\end{document}